\documentclass[11pt]{article}

\makeatletter
\newfont{\footsc}{cmcsc10 at 8truept}
\newfont{\footbf}{cmbx10 at 8truept}
\newfont{\footrm}{cmr10 at 10truept}

\makeatother
\usepackage{amsfonts,amsmath,amsthm,amssymb,mathrsfs,verbatim}

\let\n\noindent

\font\tenmsy=msbm10
\font\sevenmsy=msbm10 at 7pt
\font\fivemsy=msbm10 at 5pt
\newfam\msyfam \textfont\msyfam=\tenmsy
\scriptfont\msyfam=\sevenmsy
\scriptscriptfont\msyfam=\fivemsy

\let\s\sigma

\let\l\left
\let\r\right

\def\e{{\epsilon}}
\def\y{{\infty}}

\let\Rw\Rightarrow
\def\l{{\left}}
\def\r{{\right}}
\def\rw{\rightarrow}
\def\lrw{\leftrightarrow}

\newtheorem{theorem}{Theorem}
\newtheorem{lemma}[theorem]{Lemma}
\newtheorem{proposition}[theorem]{Proposition}

\theoremstyle{remark}

\let\s\sigma

\let\l\left
\let\r\right

\let\a\alpha
\let\b\beta

\begin{document}

\vskip18pt

 \title{\vskip60pt Multiple partitions,  lattice paths and a Burge-Bressoud-type correspondence}

\vskip18pt

\smallskip
\author{ P. Jacob and P.
Mathieu\thanks{patrick.jacob@durham.ac.uk,
pmathieu@phy.ulaval.ca.  } \\ 
\\
Department of Mathematical Sciences, \\University of Durham, Durham, DH1 3LE, UK\\
and\\
D\'epartement de physique, de g\'enie physique et d'optique,\\
Universit\'e Laval,
Qu\'ebec, Canada, G1K 7P4.
}

\vskip .2in
\bigskip

\date{August 2006}

\maketitle

\vskip0.3cm
\centerline{{\bf ABSTRACT}}
\vskip18pt

A bijection is presented between (1):  partitions with conditions  $f_j+f_{j+1}\leq k-1$ and $  f_1\leq i-1$,  where $f_j$ is the frequency of the part $j$ in the partition, 
 and (2): sets of $k-1$ ordered partitions $(n^{(1)}, n^{(2)}, \cdots ,n^{(k-1)})$ such that $n^{(j)}_\ell \geq n^{(j)}_{\ell+1} + 2j$ and $ n^{(j)}_{m_j} \geq j+ {\rm max}\, (j-i+1,0)+
2j (m_{j+1}+\cdots +  m_{k-1})$, where $m_j$ is the number of parts in $n^{(j)}$. This bijection entails an elementary and constructive proof of the Andrews multiple-sum enumerating partitions with frequency conditions.  A very natural relation between  the $k-1$ ordered partitions and  restricted paths is also presented, which reveals   our bijection to be  a modification of Bressoud's version of the Burge correspondence.

\newpage

\let\Rw\Rightarrow
\let\rw\rightarrow
\let\l\left
\let\r\right
\let\s\sigma
\let\ka\kappa
\let\de\delta

\section{Introduction}

Let $F_{k,i}(n,m)$ be the number of  partitions $(p_1,\cdots , p_m)$ of weight $n\; (=\sum p_l)$ with `difference 2 at distance $k-1$' and  at most $i-1$ occurrences of 1, namely, partitions satisfying 
\begin{equation}\label{difk}
p_l\geq p_{l+k-1}+2 \quad {\rm with} \quad p_{m-i+1}\geq 2 \quad {\rm and} \quad  p_m\geq 1.
\end{equation}
The generating function for the  numbers $F_{k,i}(n,m)$ is
 \cite{Andrr,Andr}
\begin{equation}\label{AG}
F_{k,i}(z;q)= \sum_{n,m\geq 0} F_{k,i}(n,m) q^n z^m= \sum_{m_1,\cdots,m_{k-1}=0}^\y {
q^{N_1^2+\cdots+ N_{k-1}^2+N_i+\cdots +N_{k-1} }z^{N_1+\cdots +N_{k-1} } \over (q)_{m_1}\cdots (q)_{m_{k-1}}
},
\end{equation}
with $
 N_j= m_j+\cdots
+m_{k-1}$ and $(q)_n= \prod_{i=1}^n (1-q^i)$.

An elementary parafermionic proof of  Andrews' multiple-sum  (\ref{AG}) has been presented recently  \cite{JM.A} (and found afterward to be also contained in  \cite{Geo1,Geo2} in a different setting.)
The argument of \cite{JM.A}  relies on conformal-field-theoretical methods. The idea of the proof  is as follows. For a particular  irreducible representation of the parafermionic algebra indexed by an integer $1\leq i\leq k$, we gave two bases of  states. One is expressed in terms of the modes of an operator (spanning the infinite dimensional parafermionic algebra) and the other is expressed in terms of the modes of $k-1$ distinct operators (spanning the same algebra). The elements of each basis can be put in one-to-one correspondence with partitions. 
In the first case, each state  is in correspondence with  a single  partition satisfying the condition (\ref{difk}). In the second case,  each state  is in correspondence with a set of $k-1$ ordered partitions of respective lengths $m_1,\cdots , m_{k-1}$:
\begin{equation}\label{mul}
 (n^{(1)}, n^{(2)}, \cdots ,n^{(k-1)})\qquad {\rm with}\qquad 
n^{(j)}= (n^{(j)}_1, \cdots , n^{(j)}_{m_j})\;.
\end{equation}
The parts within a given partition $n^{(j)}$ are all distinct, subject to the condition
\begin{equation}\label{difone}n^{(j)}_\ell \geq n^{(j)}_{\ell+1} + 2j\;.  
\end{equation}
There are further constraints that relate the parts within different partitions, which are
\begin{equation}\label{bon}   n^{(j)}_{m_j} \geq j+ {\rm max}\, (j-i+1,0)+
2j (m_{j+1}+\cdots +  m_{k-1})=: \Delta_{(i;j)}\; .
\end{equation}
The parafermionic aspect of the proof amounts to establish the equivalence of two bases. It implies that
 $G_{k,i}(n,m)$, enumerating the multiset (\ref{mul}) subject to  (\ref{difone})-(\ref{bon}), where
\begin{equation}
n= \sum_{j=1}^{k-1} \sum_{l=1}^{m_j} n^{(j)}_l\qquad {\rm and}\qquad   m= \sum_{j=1}^{k-1} jm_j\; ,
\end{equation}
is equal to $F_{k,i}(n,m)$. 
The next step is  to show that the generating function for the numbers $G_{k,i}(n,m)$ 
is precisely the multiple sum given in (\ref{AG}). 
This step is elementary. Indeed, the number of partitions of length  $m_j$ satisfying (\ref{difone}) is
 \begin{equation*}
 {q^{jm_j(m_j-1)}\over (q)_{m_j}} z_j^{m_j}  ,
\end{equation*}
where $z_j$ is introduced as an extra variable to keep track of the length. Taking care of condition (\ref{bon}) amounts to shifting all  $m_j$ parts of $n^{(j)}$ by the constant $\Delta_{(i;j)} $. This modifies the exponent in the numerator of the previous generating function by $m_j\Delta_{(i;j)} $: \begin{equation*}
{q^{jm_j(m_j-1)+ m_j  \Delta_{(i;j)} }\over (q)_{m_j}} z_j^{m_j}  .
\end{equation*}
By summing over all values of $m_j$ for $j=1,\ldots , k-1$ and by  setting $z_j=z^j$, one recovers  the multiple-sum expression in (\ref{AG}) \cite{JM.A}.

The essential aim of this article is to provide
a bijective proof of the equality of the numbers $F_{k,i}(n,m)$ and $G_{k,i}(n,m)$. 

 \begin{theorem}\label{F=G} Let $F_{k,i}(n,m)$ be the number of partitions of $n$ into $m$ parts satisfying (\ref{difk}) and let $G_{k,i}(n,m)$ be the number of ordered multiple partitions (\ref{mul}) satisfying (\ref{difone}) and (\ref{bon}). Then  $F_{k,i}(n,m)= G_{k,i}(n,m)$.\end{theorem}

A bijective proof of Theorem \ref{F=G}  is displayed in Section 3. The rationale  underlying this bijection is rooted in the relationship between restricted lattice paths and multiple partitions which is presented first in Section 2. That this relation is a one-to-one correspondence follows from Theorem \ref{F=G} and the Burge correspondence as interpreted in \cite{AnB,BreL}:

 \begin{proposition}\label{F=P} {\rm \cite{Bu}} Let $P_{k,i}(n,m)$ be the number of lattice paths  of weight $n$,  total charge $m$, starting at $(0,k-i)$, and with peaks of height not larger than $k-1$. Then $P_{k,i}(n,m)= F_{k,i}(n,m)$, where the numbers $F_{k,i}(n,m)$ are defined in Theorem \ref{F=G}.\end{proposition}

\n We explain the  definition of $P_{k,i}(n,m)$ in Section 2.  Note that in \cite{Bu}, the following equivalent definition of $F_{k,i}(n,m)$  is also used: 
 $F_{k,i}(n,m)$ is the number of partitions of $n$ into $m$ parts  with 
 \begin{equation}\label{fre}
f_j+f_{j+1}\leq k-1\qquad {\rm and}\qquad  f_1\leq i-1\;, 
\end{equation}
 where $f_j$ is the frequency of the part $j$ in the partition. The equivalence between the two definitions is easily verified:  if a partition satisfying (\ref{difk}) contains the sequence $(j+1,\ldots,j+1,j,\ldots j)$, the multiplicity of $j+1$ plus that of $j$ cannot be larger than $k-1$ to respect the `difference 2 at distance $k-1$' condition. Moreover, $p_{m-1+1}\geq 2$ just means that $f_1\leq i-1$.  The frequency characterization (\ref{fre}) is used explicitly in Section 3.

The relation between paths and partitions with frequency conditions is used as an intermediate step in \cite{Bu} where the objective is  to provide a simple bijective proof between these restricted partitions and partitions with successive ranks in a prescribed interval (cf. \cite{Bu} for the precise definition of this second type of partitions, which is not needed here). Note that in \cite{Bu}, the stated result is equivalent to $\sum_m P_{k,i}(n,m)= \sum_m F_{k,i}(n,m)$. But given the charge concept introduced in \cite{BreL} (under a different name), the sum over $m$ is seen to be superfluous.
That paths are interesting by themselves has been stressed in \cite{AB,AnB,BreL}. In particular, Rogers-Ramanujan-type identities are sometimes more transparent when the restrictions appropriate to partitions enumerated by the sum side of the identities are formulated  in terms of paths. That  makes relevant  the search for alternative and more intuitive links between paths and restricted partitions. 
In \cite{BreL}, Bressoud has presented such  a direct reformulation of the Burge  bijection. He has shown that a  path is essentially a sequence of peaks for which  their  height (more precisely, their charge) is related to the number of parts in sequences of the form  $(j+1,\ldots,j+1,j,\ldots j)$ within the partition. However, this reformulation  is not quite correct and the proper version is presented here (in Section 3), formulated in terms of multiple partitions.  This result is  interesting by itself.

Back to our main point: Theorem \ref{F=G} implies that the generating function for partitions with frequency conditions (\ref{fre}) is also the generating function for the multiple partitions enumerated by $G_{k,i}(n,m)$. As already stressed, obtaining the generating function for these numbers  $G_{k,i}(n,m)$ is elementary. The bottom line is thus an elementary and constructive route for obtaining $F_{k,i}(z;q)$.

\section{Lattice paths, charged clusters and multiple partitions}

The objective of this section is to show that  multiple partitions are naturally related to restricted lattice paths  \cite{BreL}. 
Restricted lattice paths  are defined in the first quadrant of   an integer square lattice as  follows.
A path starts at a non-negative  integer  position $a$ on the $y$ axis and terminates on the $x$ axis.
The possible moves are 
either  
\begin{equation*}  
\a:(i,j)\rw(i+1,{\rm max}\, (0,j-1))
\qquad {\rm 
or }\qquad
\b:(i,j)\rw(i+1,j+1).
\end{equation*} 
The restriction condition is that the height ($y$-coordinate) of the peaks  cannot be larger than $k-1$.
The weight of a path is  the sum of the $x$-coordinate of all the peaks.  The charge (called the relative height in \cite{BreL}) of 
a peak with coordinates $(i,j)$ is the largest integer $c$ such that we can find two points $(i',j-c)$ and $(i'',j-c)$ on the path  with $i'<i<i''$ and such that between these two points there are no peak of height larger than $j$ and every peak of height equal to $j$ has weight larger than $i$ \cite{BP}.
The total charge of a path is the sum of the charges of all its peaks.

A path is fully characterized by its initial
 vertical  position $k-i$  and a binary word in $\a$ and $\b$ \cite{BreL}. Equivalently, the path is  specified  by its initial point and the sequence of its peaks read from right to left, together with their respective charge. The peak specification $(x_j,c_j)$, where $x_j$ is the $x$-coordinate and $c_j$ the charge,  is conveniently written in the form $x_j^{(c_j)}$. For instance, with $i=k\geq 4$, we have 
\begin{equation*} \a\b^2\a^2\b\a^2\b^2\a\b\a^2\b\a\b^3\a^3 \lrw \, 19^{(3)}\, 15^{(1)}\, 12^{(1)} \, 10^{(2)} \, 6^{(1)} \, 3^{(2)}\;.  \end{equation*} 
From now on, it will be understood that the ordering matters when specifying a sequence of $x_j^{(c_j)}$.

 The basic  characteristics of a path are captured by the following lemma, whose proof is immediate given the graphical description of the corresponding situations.

\begin{lemma} \label{path} For a restricted path with original vertical position $k-i$, the following two conditions must be satisfied:

\n (1) A peak of charge $j$ has minimal $x$-coordinate $ j+ 
{\rm max}\, (j-i+1,0)$. 

\n (2) If between two  peaks $ x^{(i)}$ and $ {y}^{(j)}$ there are peaks all with charge lower than ${\rm min}\, (i,j)$ and whose total charge sums to $c$, then
\begin{equation}\label{dist} x- y \geq r_{ij}+\chi_{i>j}+2c\;, \end{equation} 
where 
\begin{equation} \label{rij}
 r_{ij} := 2\, {\rm min} \; (i,j)\;, 
 \end{equation}
 and $\chi_{b} =1$ if $b$ is true and 0 otherwise. The special case $c=0$ describes the minimal separation between two adjacent peaks.
 \end{lemma}

\n For instance,  the difference in weights between each adjacent pair of peaks of the  sequence  $21^{(3)}\, 19^{(1)}\, 16^{(2)}\, 12^{(3)}$ satisfies (\ref{dist}). However, the difference between the two extremal peaks should be at least $6+2(1+2)=12$, meaning that this sequence of peaks does not represent a path (but $24^{(3)}\, 19^{(1)}\, 16^{(2)}\, 12^{(3)}$ would).

We now introduce a formal operation that describes the interchange of two adjacent peaks $ x^{(i)} {y}^{(j)}$.  It is defined as follows:
\begin{equation} \label{com}
 x^{(i)} {y}^{(j)} \quad \rw \quad ({y}+r_{ij}) ^{(j)} (x-r_{ij})^{(i)} \;, 
 \end{equation}
 where $r_{ij}$ is defined in (\ref{rij}).
This  operation preserves the individual values of the charge and also the sum of the weights.

It is important to stress that after a sequence of interchanges, such that $\{ x_j^{(c_j)} \}\rw
 \{ {x'_j}^{(c_j)} \}$, the new values of $x_j'$ are no longer necessarily decreasing and they no longer  correspond to peak positions in a modified path.  A specific terminology for the resulting numbers $ {x'_j}^{(c_j)} $ is thus required. We will refer to ${x'_j}^{(c_j)}$ as a cluster of charge $c_j$ and  weight $x'_j$.
 Clearly, the interchange operation is defined for any two clusters and it is not restricted to peaks.

Let $P_{k,i}(n,m)$ and  $G_{k,i}(n,m)$ be  defined as in Proposition \ref{F=P} and Theorem \ref{F=G} respectively. The rest of this section is devoted to establishing the following proposition:

 \begin{proposition}\label{P=G} We have $P_{k,i}(n,m)= G_{k,i}(n,m)$.\end{proposition}

\n {\it Remark}: Let us point out at first that Theorem \ref{F=G} and Proposition \ref{F=P} are independent of the relationship between paths and multiple which is  presented below. And taken together, they entail Proposition \ref{P=G}. For this reason, 
we confine ourself to the presentation of a  sketch of a direct proof.

\begin{proof} 
Let us first present the correspondence between a path and a multiple partition $(n^{(1)},\cdots,n^{(k-1)})$. Start from a path expressed as an ordered sequence of clusters and interchange the clusters using (\ref{com}) until they   are  all ordered with increasing charge (from left to right) and, within each sequence of clusters with identical charge, with decreasing weight. The weights of the clusters of  charge $j$ then form the parts of the partition $n^{(j)}$. We view the multiple partition as the  canonical rewriting of the original path.
 For instance, we have
\begin{equation*}  19^{(3)}\, 15^{(1)}\, 12^{(1)} \, 10^{(2)} \, 6^{(1)} \, 3^{(2)} \rw 
 17^{(1)}\, 14^{(1)}\, 10^{(1)} \, 12^{(2)} \, 7^{(2)} \, 5^{(3)} \,.\end{equation*} 
The multiple partition thus obtained is
$  n^{(1)} = (17,14,10)$, $n^{(2)} = (12,7)$ and $n^{(3)} = (5).$

It is easy to see that the correspondence between  a path and a multiple partition is a well-defined map and that the conditions (\ref{difone}) and (\ref{bon}) are satisfied by construction, being 
 ensured by the definition of the interchange operation and the basic properties of a path given in Lemma \ref{path}.

The inverse operation amounts to rewriting $(n^{(1)},\cdots,n^{(k-1)})$ as a sequence of clusters and then reordering the clusters using (\ref{com})   to ensure that the conditions (\ref{dist}) are  everywhere satisfied.  For this, one first reorders the clusters such that the weight are decreasing from left to right. 
Then, for every pair of adjacent clusters $x^{(i)} y^{(j)}$ that do not satisfy (\ref{dist}), one performs an interchange. The condition (\ref{dist}) is necessarily satisfied after the interchange:
\begin{equation*}x- y <r_{ij}+\chi_{i>j} \quad \Rw\quad  (y+r_{ij})-(x-r_{ij} ) \geq r_{ij}-\chi_{i>j} +1\geq r_{ij}+ \chi_{j>i} \;. \end{equation*} 
 We thus perform all the required  interchanges of  adjacent pairs by selecting always the leftmost problematic pair. Once all adjacent pairs are transformed in order to satisfy (\ref{dist}), we then look for pairs of clusters of charge $j\geq 2$ that are separated by clusters of lower charges and reorder those pairs that do not satisfy (\ref{dist}). This operation might produce violations of (\ref{dist}) among adjacent pairs and the process of reordering the latter is then done anew. The procedure is completed once (\ref{dist}) is satisfied for all pairs of clusters.  This process terminates. The resulting configuration is the corresponding path. 
\end{proof}

\n {\it Remark}:  We note a simple mechanical characterization of a path in terms of an equilibrium condition (an energy minimum) for a systems of  interacting clusters. Let us thus  define an  interaction energy $E$ for a sequence of clusters by summing over all pairwise  interaction energies $\e_{ij}$ defined as follows: 
For two clusters $x^{(i)}$ and $y^{(j)}$ with $x>y$ separated by a sequence of intermediate peaks each with charge lower than min $(i,j)$ and whose total charge sums to $c$, this energy $\e_{ij}$ is 
 \begin{equation*}
 \e_{ij}={\rm max}\, \big(r_{ij}+\chi_{i>j}+2c-(x-y), 0\big)\,.
\end{equation*}
 If  in between $x^{(i)}$ and $y^{(j)}$ there is at least one  cluster of charge $\geq$ min $(i,j)$, we set $\e_{ij}=0$.
The interaction energy of a sequence of clusters is  then defined to be
 $E= \sum_{i < j } \e_{ij}$. 
A path corresponds to  a sequence of clusters ordered by decreasing weight which  has vanishing interaction energy. The last condition is a consequence of (\ref{dist}).  Actually, the above proposition implies that, given a multiple partition, the corresponding zero-energy configuration is unique.

\n {\it Remark:}  
We stress that listing  multiple partitions is easier than listing paths.
For instance, listing all multiple partitions with $n=17$ and $m_1=1,\, m_2=2$ (so that $m=5$), with $i=k=3$, amounts to list triplets of integers $(n_1^{(1)}; n_1^{(2)}, n_2^{(2)})$ satisfying  
 \begin{equation*}
n_1^{(1)}+ n_1^{(2)}+ n_2^{(2)}=17 \quad {\rm with }\quad n_1^{(1)} \geq 5, \quad n_1^{(2)} \geq n_2^{(2)}+4,\quad {\rm and}\quad  n_2^{(2)}\geq 2 .
 \end{equation*}

\section{A Burge-Bressoud-type correspondence}

In the previous section, we have  motivated the representation of a multiple partition as an ordered sequence of clusters (where $n_\ell^{(j)}$ is called a cluster of weight $n_\ell$ and charge $j$) and  introduced a `commutation rule' for clusters, which is the interchange operation given in (\ref{com}) with $r_{ij}$ defined by (\ref{rij}). These are the only ingredients needed from Section 2 for  the description of the bijection  between  partitions with frequency conditions and  multiple partitions that is constructed here.

We now turn to the proof of Theorem \ref{F=G}.

\begin{proof} 
Let us first show how to associate a multiple partition to a partition $(p_1,\cdots , p_m)$ satisfying (\ref{difk}). Identify the sequences of $k-1$ adjacent parts such that the first part and last part of each  sequence differ by at most 1. In other words, find sequences $(p_j,\cdots, p_{j+k-2})$ such that $p_j-p_{j+k-2}\leq 1$. The condition (\ref{difk}) ensures that two sequences of this type cannot overlap. Each such sequence is then replaced by a cluster of charge $k-1$ whose weight is given by the sum of  its  parts, i.e.,
\begin{equation*}( p_j,\cdots, p_{j+k-2} ) \rw 
\left(\sum_{r=j}^{j+k-2} p_r\right )^{(k-1)}.
\end{equation*} 
Once all clusters of charge $k-1$ are constructed, we move them (preserving their ordering) to the right of the sequence formed by the remaining parts. This displacement is done using the interchange operation (\ref{com}),  by treating all parts which are crossed as clusters of charge 1. Once this is completed, one is  left with  a smaller partition and a sequence of ordered clusters of charge $k-1$ at its right. For the resulting partition, one repeats the previous analysis but with $k-1$ replaced by $k-2$. Once all clusters of charge $k-2$ are identified, they are moved to the right of the partition.
This procedure is repeated for lower-charge clusters until all clusters of charge 2 are formed and moved to the left extremity of the sequence of ordered clusters of charge $3,\cdots, k-1$. The remaining parts of the partition are the clusters of charge 1. The result is a multiple partition of the form (\ref{mul}), where the parts of $n^{(j)}$ are the weights of the  clusters of charge $j$.  
The conditions (\ref{difone})--(\ref{bon}) are immediate consequences of the construction.

Whenever $f_s+f_{s-1}=f_{s-1}+f_{s-2}$, with $f_{s-1}>0$, that is, for partitions of the form
 \begin{equation*}
(\ldots, s', \underbrace{s,\cdots , s}_{r}, \underbrace{s-1, \cdots, s-1}_{j-r}, \underbrace{s-2, \cdots, s-2}_{r}, s'', \ldots)\;,
 \end{equation*}
there is a potential ambiguity in the regrouping of parts. A cluster of charge $j$ can be obtained  
 by regrouping the sequences of $s-1$ and $s-2$ as
  \begin{equation*}
(\ldots, s', s,\cdots , s, \, (j(s-2)+j-r)^{(j)}, s'', \ldots)\;, 
\end{equation*}
or either  by regrouping the parts $s$ and $s-1$, as  
 \begin{equation*}
(\ldots, s', \left(j(s-1)+r\right)^{(j)} ,s-2, \cdots, s-2, s'', \ldots) \;.\end{equation*}
But in the latter case, by commuting the cluster through the $r$ parts equal to $s-2$, one recovers the former expression. The ambiguity in the regrouping process is thus superficial (in that the different clusterings are related by interchange) and does not affect the final multiple partition.

Let us illustrate the clustering procedure for a partition that satisfies (\ref{difk}) for $k=5$:
\begin{equation*}
(8,8,7,7,5,3,3,2,2,1,1) \rw (30^{(4)}, 5, 10^{(4)}, 1,1) \rw (7, 5,5) \,24^{(4)}\,6^{(4)}\rw 7^{(1)}\,  10^{(2)}\,24^{(4)}\,6^{(4)}.\end{equation*}
The corresponding multiple partition is thus  $ n^{(1)}= (7), n^{(2)}= (10)$ and $ n^{(4)}= (24,6)$.

The inverse operation is formulated as follows. Re-express the multiple partition  as the partition $(n_1^{(1)}\cdots  n_{m_1}^{(1)})$ followed by the sequence of clusters 
 $n_1^{(2)} \cdots n_\ell^{(j)} \cdots n_{m_{k-1}}^{(k-1)}$. Each cluster is then inserted successively  (starting with $n_1^{(2)}$ up to $n_{m_{k-1}}^{(k-1)}$) within the partition, using the interchange rule  (\ref{com}), treating again each part as a cluster of charge 1. Once inserted within the partition (at a position to be determined below), a cluster is unfolded into  the number of parts given by its charge, with parts as equal as possible. To state this precisely, consider for definitiveness the insertion of ${n}_\ell^{(j)}$. Its displacement within the partition modifies its weight to $n'_\ell$ ($n'_\ell=n_\ell  \; +$ twice  the number of interchanges performed). The cluster ${n'}_\ell^{(j)}$ is  then broken apart  into $j$ parts differing at most by 1 and whose sum is $n'_\ell$. This decomposition is unique since given $n_\ell$ and $j$, there are unique non-negative integers $s$ and $r$ such that $n'_\ell= sj+r$ with $r$ smaller than $j$; in the decomposition of $n'_\ell$, there are then $r$ parts equal to $s+1$ and $j-r$ parts equal to $s$:
 \begin{equation*}
 {n'}_\ell^{(j)} \rw (\underbrace{s+1,\cdots , s+1}_{r}, \underbrace{s, \cdots, s}_{j-r})\;.
 \end{equation*}
The position at which the cluster ${n'}_\ell^{(j)}$ is moved within the partition is determined by two criteria: 

\begin{enumerate}
\item The new sequence of numbers  that results from unfolding ${n'}_\ell^{(j)}$ must be a partition.
\item The frequency condition $f_t+f_{t-1}\leq j$ must be satisfied for all parts $t$ of this partition. Equivalently, the resulting partition must satisfy
 \begin{equation}\label{disj}
 p_l-p_{l+j}\geq 2\;.
 \end{equation}
 \end{enumerate}

Obviously, because $j\leq k-1$, (\ref{disj})  ensures the validity of the condition (\ref{difk}) at every intermediate stage of the  construction.

To show that these conditions fix the position where the cluster must be placed and unfolded, it suffices to consider the case of two clusters $x^{(i)} \, y^{(j)}$. Let us  first consider the situation where $i<j$ and express the weights as $x= is+r$ and $y = js'+r'.$ Unfolding the first cluster yields:
\begin{equation*}
 (is+r)^{(i)} \, (js'+r')^{(j)} \rw (\underbrace{s+1,\cdots , s+1}_{r}, \underbrace{s, \cdots, s}_{i-r}) \, (js'+r')^{(j)}\;. 
 \end{equation*}
 If $s>s'+2$, the second cluster is unfolded at the right of the partition  resulting into 
 \begin{equation*}
(\underbrace{s+1,\cdots , s+1}_{r}, \underbrace{s, \cdots, s}_{i-r}, \underbrace{s'+1,\cdots , s'+1}_{r'}, \underbrace{s', \cdots, s'}_{j-r'}) \;.
 \end{equation*}
This is indeed a partition  and since $s-(s'+1)>1$, (\ref{disj}) is satisfied.

If $s\leq s'$, the second cluster must be commuted at the beginning of the partition and then unfolded. This corresponds to a case where the two clusters are interchanged, 
\begin{equation*}
(is+r)^{(i)} \, (js'+r')^{(j)} \rw (js'+r'+2i)^{(j)} \, (i(s-2)+r)^{(i)} \;,
 \end{equation*} and then unfolded. The condition $s\leq s'$ ensures that the resulting sequence is non-increasing. If $s<s'$, the condition (\ref{disj}) is automatically verified.  The case $s=s'$ needs a slightly more careful inspection. If $r'+2i<j$, the unfolding yields:
\begin{equation*}
(\underbrace{s+1,\cdots , s+1}_{r'+2i}, \underbrace{s, \cdots, s}_{j-r'-2i}, \underbrace{s-1,\cdots , s-1}_{r}, \underbrace{s-2, \cdots, s-2}_{i-r}) \;.
 \end{equation*}
 There is a potential problem with the length of the sequences of $s$ and $s-1$: $f_s+f_{s-1}=j+r-r'-2i$,  which is greater than $j$ if $r>r'+2i$. But since $r<i$, the inequality $r>r'+2i$ can never be satisfied. Finally, when $r'+2i\geq j$, (\ref{disj}) is directly verified.

The remaining cases, namely $s'=s-1$ or $s-2$, require a more detailed analysis. We will spell out the details  pertaining to  $s'=s-1$ in order  to illustrate the procedure and just state the results for $s'=s-2$.

By unfolding the cluster $(j(s-1)+r')^{(j)}$ at the end of the partition, we get:
 \begin{equation*}
(\underbrace{s+1,\cdots , s+1}_{r}, \underbrace{s, \cdots, s}_{i-r}, \underbrace{s,\cdots , s}_{r'}, \underbrace{s-1, \cdots, s-1}_{j-r'}) \;.
 \end{equation*}
 This is indeed a non-increasing sequence. However, because $f_s+f_{s-1}=j+i-r>j$ (since $i-r>0$), the second criterion is violated. Suppose instead  that we interchange the two clusters, getting 
 $ (j(s-1)+r'+2i)^{(j)} \, (i(s-2)+r)^{(i)} $,
 and unfold them. We then need to consider three distinct situations according to the value of $r'+2i$:
 
 \n i) $r'+2i<j$. In that case, the unfolding reads
  \begin{equation*}
 (\underbrace{s,\cdots , s}_{r'+2i}, \underbrace{s-1, \cdots, s-1}_{j-r'-2i}, \underbrace{s-1,\cdots , s-1}_{r}, \underbrace{s-2, \cdots, s-2}_{i-r}) \;.
 \end{equation*}
 Here one sees that $f_s+f_{s-1}=j+r$ which is $>j$ when  $r>0$,  meaning that this unfolding position is not the appropriate one.
On the other hand, for $r=0$, the two criteria are verified.

\n ii) $j\leq r'+2i<2j$: Setting  $r''= r'+2i-j$, we have
   \begin{equation*}
  (\underbrace{s+1,\cdots , s+1}_{r''}, \underbrace{s, \cdots, s}_{j-r''}, \underbrace{s-1,\cdots , s-1}_{r}, \underbrace{s-2, \cdots, s-2}_{i-r}) \;.
 \end{equation*}
Here  $f_s+f_{s-1}=j+r-r''>j$ whenever $r>r''$. Otherwise, the two criteria are satisfied.

\n iii) $2j\leq r'+2i<3j$: With $r'''=r'+2i-2j$, one obtains
   \begin{equation*}
 (\underbrace{s+2,\cdots , s+2}_{r'''}, \underbrace{s+1, \cdots, s+1}_{j-r'''}, \underbrace{s-1,\cdots , s-1}_{r}, \underbrace{s-2, \cdots, s-2}_{i-r}) \;.
 \end{equation*}
This partition satisfies (\ref{disj}).

Summarizing, when $s'=s-1$, unfolding the charge $j$ cluster at the end of the partition never works; if 
  $r\leq {\rm max}\; (0,r'+2i-j)$
 the unfolding must be done at the beginning of the partition. However, when $r> {\rm max}\; (0,r'+2i-j)$, none of these two positions for unfolding the charge $j$ cluster  satisfies the above conditions. The charge $j$ cluster must then be placed {\it within} the partition of $i$ parts before being unfolded.

Note that for the problematic cases,  $r$ cannot be zero. This means that the partition into $i$ parts has some parts equal to $s+1$ and some equal to $s$.  It is very simple to verify that  unfolding the charge $j$ cluster at any  internal position other than the $(r+1)$-th one --  which is precisely the point where the parts of the original partition differ by 1 -- violates (\ref{disj}).  It remains to verify that the final  possibility, that is,  placing the charge $j$ cluster directly after the sequence of $r$ parts equal to $s+1$ and unfolding it, produces a partition with the proper frequency condition.
This results in 
\begin{equation*}
(\underbrace{s+1,\cdots , s+1}_{r}, \underbrace{s, \cdots, s}_{r'+2i-2r}, \underbrace{s-1,\cdots , s-1}_{j-r'-2i+2r}, \underbrace{s-2, \cdots, s-2}_{i-r}) \;. 
 \end{equation*}
 One then only needs to check that the length of  the sequences of $s+1$ and $s$ and that of the sequences of  $s$ and $s-1$, are properly bounded. With $r> {\rm max}\; (0,r'+2i-j)$, it follows that 
 \begin{equation*} f_{s+1}+f_{s}=r'+2i-r<r'+2i+{\rm min}\; (0,j-r'-2i)= {\rm min}\; (r'+2i,j)\leq j\;. 
  \end{equation*}
  Similarly, because $r<i$, one has
  \begin{equation*} f_{s-1}+f_{s-2}=j-r'-(i-r)<j\;,
  \end{equation*}
  as it should.

For $s'=s-2$, a similar analysis shows that the unfolding of $y^{(j)}$ is done at the end of the partition when  $i-r+r'\leq j$ and at the beginning if   $i-r+r'>j$ and $r\leq {\rm max}\; (0,r'+2i-2j)$. Otherwise, the unfolding is done after placing $y^{(j)}$ directly at the right of  the sequence of the $r$ consecutive $s+1$ of the partition.

 Finally, consider the situation where $i=j$, in which case,  $x-y \geq 2j$.  With $x= js+r$ and $y = js'+r'$, this requires either  $s>s'+2$ or $s=s'+2$ and $r\geq r'$. The direct unfolding of the two clusters yields
 \begin{equation*}
(\underbrace{s+1,\cdots , s+1}_{r}, \underbrace{s, \cdots, s}_{j-r}, \underbrace{s'+1,\cdots , s'+1}_{r'}, \underbrace{s', \cdots, s'}_{j-r'}) \;.
 \end{equation*}
 The condition (\ref{disj}) is satisfied. 
 
 Note that the superficial ambiguity mentioned previously in the regrouping of clusters has its counterpart in the inverse operation: there might be two positions at which a cluster could be placed and unfolded. But whenever this is the case, the resulting two partitions are identical.
\end{proof}

The simplest example that shows that the mere reordering and unfolding  of clusters does not always produce a partition satisfying (\ref{difk}) is
$7^{(2)}\, 8^{(4)}$, for $k=5$. Unfolding directly the two clusters yields $(4,3,2,2,2,2)$ which violates (\ref{difk}). Interchanging the clusters and unfolding them yields $12^{(4)}\, 3^{(2)}\rw (3,3,3,3,2,1)$, which  suffers from the same deficiency.  So we have to unfold the charge 4 cluster within the partition obtained by unfolding the charge 2 cluster: $(4,3)\, 8^{(4)}\rw (4, 10^{(4)}, 1) \rw (4,3,3,2,2,1)$.

\section{Concluding remarks: links  to other works}

In \cite{BreL} Bressoud  provided a reformulation of the Burge correspondence  between lattice paths 
and  partitions with frequency conditions that reveals its  essence 
as a kind of blowing up of each peak of weight $x$ and charge $j$ into $j$ parts differing  at most by 1 and whose sum is $x$. Our bijection, formulated in terms of multiple partitions instead of lattice paths, is thus a sort of modified version  of  Bressoud's correspondence. But in addition to being a variant, our correspondence sharpens  that of Bressoud.  
The need for a  rule that goes beyond the sole interchange of the original clusters before their unfolding  
(illustrated by the last example of the previous section  that corresponds to  the path $12^{(4)}\, 3^{(2)}$) supports this point.

In that regard, let us indicate that our key  interchange operation  (\ref{com}) can be recognized as the shuffle operation introduced in \cite{BreL} (Section 4), but without imposing the defining application criterion presented there.  From our perspective, however,  this operation is a simplified abstraction  of a commutation relation  involving  quantum operators used in \cite{JM.A}.

For given lengths $m_1,\cdots , m_{k-1}$, there is a trivial bijection between the ordered partitions 
$(n^{(1)},\cdots,n^{(k-1)}) $ of weight $n$ and of representations of $n$ in the form $ N_1^2+\cdots+ N_{k-1}^2+N_i+\cdots +N_{k-1}$  plus $k-1$ partitions with at most $m_1, m_2,\cdots , m_{k-1}$ parts. (The previous  number  is precisely the minimal weight of the set of ordered partitions with the specified length content.)
Let $D_{k,i}(n)$ be the number of representations of such  $n$. Burge \cite{Bur} has provided a bijection between elements enumerated by $D_{k,i}(n)$ and $F_{k,i}(n)=\sum_m F_{k,i}(n,m)$ (this is his third-way correspondence). In that perspective, our result can be regarded as a new derivation of the Burge's result (in a sharpened version in that the $m$ dependence is taken into account).

Finally, we indicate that the equivalence between the generating function of 
lattice paths and that of restricted partitions can be extended to restricted jagged partitions and lattice paths with special  conditions \cite{JPath}. (Jagged partitions are presented in \cite{FJM.R} and their restricted versions are considered in \cite{FJM.E}. These partitions have appeared originally in a physical context \cite{JM}). The underlying bijection between jagged partitions and paths, formulated in terms of multiple partitions, is worked out in Appendix A of \cite{Jstat}. Finally, Theorem \ref{F=G} is a key ingredient in  the proof of equivalence presented in \cite{Path}Ê between the  lattice paths considered here (which could rightly be called the Bressoud lattice paths) and another class of paths which are  naturally defined from the statistical models introduced in \cite{ABF} Ê(the so-called RSOS paths). 
\vskip0.3cm
\noindent {\bf ACKNOWLEDGMENTS}

We thank D. Ridout for a critical reading of the manuscript and clear-sighted comments. The work of PJ is supported by EPSRC and partially  by the EC
network EUCLID (contract number HPRN-CT-2002-00325), while that of  PM is supported  by NSERC.

\end{document}